\documentclass{article}
\usepackage[top=30truemm,bottom=30truemm,left=25truemm,right=25truemm]{geometry}
\usepackage{amsmath,amssymb,amsfonts,mathtools,accents,bm,mathrsfs,bibentry,lscape}
\usepackage{enumerate}
\usepackage{tabu}
\newtheorem{theorem}{Theorem}

\newtheorem{cor}{Corollary}
\newtheorem{assumption}{Assumption}

\newtheorem{proposition}{Proposition}

\title{Asymptotic properties of parametric and nonparametric probability density estimators of sample maximum}
\author{
Taku MORIYAMA\\
Department of Management of Social Systems and Civil Engineering, Tottori University}
\date{}

\begin{document}
\maketitle

\begin{abstract}
Asymptotic properties of three estimators of probability density function of sample maximum $f_{(m)}:=mfF^{m-1}$ are derived, where $m$ is a function of sample size $n$. One of the estimators is the parametrically fitted by the approximating generalized extreme value density function. However, the parametric fitting is misspecified in finite $m$ cases. The misspecification comes from mainly the following two: the difference $m$ and the selected block size $k$, and the poor approximation $f_{(m)}$ to the generalized extreme value density which depends on the magnitude of $m$ and the extreme index $\gamma$. The convergence rate of the approximation gets slower as $\gamma$ tends to zero. As alternatives two nonparametric density estimators are proposed which are free from the misspecification. The first is a plug-in type of kernel density estimator and the second is a block-maxima-based kernel density estimator. Theoretical study clarifies the asymptotic convergence rate of the plug-in type estimator is faster than the block-maxima-based estimator when $\gamma> -1$. A numerical comparative study on the bandwidth selection shows the performances of a plug-in approach and cross-validation approach depend on $\gamma$ and are totally comparable. Numerical study demonstrates that the plug-in nonparametric estimator with the estimated bandwidth by either approach overtakes the parametrically fitting estimator especially for distributions with $\gamma$ close to zero as $m$ gets large.
\end{abstract}

\section{Parametric density estimation of sample maximum}
Let $X_1,X_2,\cdots$ be independent and identically distributed random variables with a continuous distribution function $F$, which belongs to either one of the Hall class, the Weibull class and the bounded class shown in the following assumption. 
\begin{assumption}
It holds that either one of the followings {\rm (i)} $\exists \alpha>0, \exists \beta>0, \exists A>0, \exists B\neq0$ \text{s.t.} 
$$x^{\alpha+\beta} \{1 - F(x) - A x^{-\alpha} (1 + B x^{-\beta} ) \} \to 0~~~ \text{as}  ~~~ x \to \infty,$$
{\rm (ii)} $\exists \kappa>0, \exists C>0$ \text{s.t.}
$$\exp(Cx^{\kappa})\{1 - F(x) - \exp(-Cx^{\kappa})\} \to 0~~~ \text{as}  ~~~ x \to \infty,$$
{\rm (iii)} $\exists\mu<-2, \sigma<0, \exists D\neq0, \exists E\neq0, \exists x^*\in \mathbb{R}$ s.t. 
$$(x^* -x)^{\mu(1-\sigma)} \{1 -F(x) - (x^* -x)^{-\mu} (D + E (x^* -t)^{\mu\sigma})\} \to 0  ~~~ \text{as}  ~~~ x \uparrow x^*.$$
\end{assumption}
The Pareto distribution, $t$-distribution, Burr distribution, and extreme value distribution with $\gamma>0$ belong to the Hall class. The Hall class of distributions are heavily tailed. Methods for checking whether the data are heavily tailed are provided in Resnick [1]. The gamma, logistic, and log-normal distributions belong to the Weibull class of distributions. The uniform distributions, beta distributions and the reversed Burr distribution belong to the bounded class of distributions. The three class of distributions are widely used as representatives of the distributions belonging to the  domain of attractions of the extreme value distribution (see Beirlant et al. [2]; Segers [3]).

Under Assumption 1, it follows from the Fisher-Tippett-Gnedenko theorem that
\begin{align}
F^m(x) - G_{\gamma}\left(\frac{x -b_m}{a_m}\right) \to 0 ~~~\text{as} ~~~ m \to \infty,
\end{align}
where $G_{\gamma}$ is the generalized extreme value distribution, 
$$\gamma := 
\begin{dcases}
{\alpha}^{-1} ~~~ &\text{for} ~~~ {\rm (i)}\\
0 ~~~ &\text{for} ~~~ {\rm (ii)}\\
{\mu}^{-1} ~~~ &\text{for} ~~~ {\rm (iii)},
\end{dcases}$$
\begin{align*}
a_m := 
\begin{dcases}
\gamma (Am)^{\gamma} ~~~ &\text{for} ~~~ {\rm (i)}\\
\kappa^{-1} C^{-1/\kappa} (\ln m)^{-\theta} &\text{for} ~~~ {\rm (ii)} \\
-\gamma (Dm)^{\gamma} ~~~ &\text{for} ~~~ {\rm (iii)},
\end{dcases} ~~~ 
b_m := 
\begin{dcases}
(Am)^{\gamma} ~~~ &\text{for} ~~~ {\rm (i)}\\
(C^{-1} {\ln m})^{1/\kappa} ~~~ &\text{for} ~~~ {\rm (ii)}\\
x^* -(Dm)^{\gamma} ~~~ &\text{for} ~~~ {\rm (iii)}.
\end{dcases}
\end{align*}
and $\theta :=1-(1/\kappa)$. Moriyama [4] discussed the distribution estimation of the sample maximum $\max \{X_{n+1}, \cdots, X_{n+m}\}$, i.e, $F^m$. Let the continuous density $f$ exist. Then, the sample maximum density ({\it smd}) is given by $f_{(m)} := m f F^{m-1}$. We suppose $m=m_n \to \infty$ as $n \to \infty$ and  investigate asymptotic properties of several {\it smd} estimators. Beirlant and Devroye [5] discussed non-zero lower bounds of the total variation and the supremum norm in {\it smd} estimation for general settings. Novak [6] proved that there are no consistent distribution estimators of sample maximum unless $m/n=o(1)$ without any assumptions. 

Since the convergence (1) is uniform (see Beirlant et al. [2]), the following convergence on {\it smd} holds
\begin{align*}
f_{(m)}(x) - g_{\bm{\gamma}_m}(x) \to 0 ~~~ \text{as} ~~~ m \to \infty,
\end{align*}
where
\begin{align*}
g_{\bm{\gamma}_m}(x) :=& \frac{1}{a_m} g_{\gamma}\left(\frac{x -b_m}{a_m}\right) \\
g_{\gamma}(x) =&
\begin{dcases}
\{1+\gamma x\}^{-(1/\gamma)-1} \exp(-\{1+\gamma x\}^{-1/\gamma}) , ~~~ &1+\gamma x>0, ~~~ \text{for} ~~~ {\rm (i), (iii)}\\
\exp(-x) \exp(-\exp(-x)), ~~~ &x \in \mathbb{R}, ~~~ \text{for} ~~~ {\rm (ii)}.
\end{dcases}
\end{align*}
We estimate the parameters $\bm{\gamma}_m :=(\gamma, a_m, b_m)$ by the maximum likelihood estimation based on the block maxima method (BM), where the block size is $k$ and ${}^{\exists}N \in \mathbb{N}$ s.t. $n=N\times k$. $k$ can be different from $m$; however, we do not consider block size selection. 
\begin{assumption}
Under Assumption 1, either $(M_n \vee K_n) \to 0$, $(M_n \wedge K_n) \to \infty$, or ${}^{\exists}\delta>0$ s.t. $M_n \to \delta$ and $K_n \to \delta$ holds, where 
\begin{align*}
M_n := 
\begin{dcases}
A m x_{n}^{-\alpha} ~ &\text{for} ~ {\rm (i)}\\
m \exp(-C x_n^{\kappa}) ~ &\text{for} ~ {\rm (ii)}\\
D m (x^* -x_{n})^{-\mu} ~ &\text{for} ~ {\rm (iii)},
\end{dcases} ~~~
K_n :=
\begin{dcases}
A k x_{n}^{-\alpha} ~ &\text{for} ~ {\rm (i)}\\
k^{\kappa} \exp (-\kappa C^{1/\kappa} (\ln k)^{\theta} x_n ) ~ &\text{for} ~ {\rm (ii)}\\
D k (x^* -x_{n})^{-\mu} ~ &\text{for} ~ {\rm (iii)}.
\end{dcases}
\end{align*}
\end{assumption}
\begin{proposition}
Under Assumptions 1 and 2,
$$\widetilde{\tau}_n := f_{(m)}(x_n) - g_{{\bm{\gamma}}_k}(x_n) \to 0.$$
\end{proposition}

\begin{assumption}
Under Assumption 1, ${}^{\exists}\lambda \in \mathbb{R}$ s.t. $\lambda_n \to \lambda$, where 
$$\lambda_{n} := 
\begin{dcases}
k m^{-2\beta} ~~~ &\text{for} ~~~ {\rm (i)}\\
k (\ln m)^{-2} ~~~ &\text{for} ~~~ {\rm (ii)}\\
k m^{2\sigma} ~~~ &\text{for} ~~~ {\rm (iii)}.
\end{dcases}$$
\end{assumption}

\begin{assumption}$N^{-1/2} k^{-\gamma} x_n \to 0$.
\end{assumption}

Set the MLE $\widehat{\bm{\gamma}}_k$. Regarding the asymptotic normality of the parametric estimator (PE) $g_{\widehat{\bm{\gamma}}_k}(x_n)$, the following theorem holds. 
\begin{theorem}
Under Assumptions 1,2, 4 and 5,
$$\sqrt{N} k^{\gamma} (\widetilde{\bm{\eta}}_n^{\mathsf{T}} I_{0}^{-1} \widetilde{\bm{\eta}}_n)^{-1/2} \left\{(f_{(m)}(x_n) - g_{\widehat{\bm{\gamma}}_k}(x_n)) - (\widetilde{\tau}_n + N^{-1/2} k^{-\gamma} \lambda_n \widetilde{\bm{\eta}}_n^{\mathsf{T}} I_{0}^{-1} \bm{b})\right\}$$ 
converges in distribution to the standard normal, where $\widetilde{\bm{\eta}}_n := (\widetilde{s}_n, \widetilde{t}_n, \widetilde{u}_n)^{\mathsf{T}}$, 
\begin{align*}
\widetilde{s}_n &:= K_n^{1+\gamma}{\exp}(-K_n) \left\{(1 -K_n)(1-K_n^{\gamma} +\gamma \ln K_n) + \gamma \left(1-K_n^{\gamma} \right) \right\}, \\
\widetilde{t}_n &:= K_n^{1+\gamma} {\exp}(-K_n) (K_n -1) \left(\frac{K_n^{\gamma} -1}{\gamma}\right), \\
\widetilde{u}_n &:=K_n^{1+\gamma} {\exp}(-K_n) K_n^{\gamma} (1+\gamma -K_n).
\end{align*} 
 $\bm{b}$ and the Fisher information matrix $I_{0}$ are given in Dombry and Ferreira [7].
\end{theorem}
This paper employs the symbols used in Moriyama [4] for clarity. We add a tilde on the symbols that are used in the same manner as in Moriyama [4] but are in fact different from those in Moriyama [4]. 
The following corollary on the convergence rate follows from Theorem 6. Proposition 3 and Theorem 10 are  proved in the Appendices, and Corollary 7 is a direct consequence of Theorem 6.
\begin{cor}
Under the assumptions of Theorem 6, $(f_{(m)}(x_n) - g_{\widehat{\bm{\gamma}}_k}(x_n))$ converges with the rate $(N^{-1/2} k^{-\gamma} \lambda_n \widetilde{\bm{\eta}}_n^{\mathsf{T}} \bm{1} + \widetilde{\tau}_n + \widetilde{\zeta}_n)$, where 
\begin{align*}
\widetilde{\zeta}_{n}:=& N^{-1/2} k^{-\gamma}\times
\begin{dcases}
 K_n^{1+\gamma} \ln K_n ~~~ &\text{for} ~~~ (M_n \vee K_n) \to 0 \\
 1~~~ &\text{for} ~~~ M_n \to \delta, ~ K_n \to \delta \\
K_n^{2(1+\gamma)} \exp(-K_n) ~~~ &\text{for} ~~~ (M_n \wedge K_n) \to \infty.
\end{dcases}
\end{align*}
\end{cor}

The misspecification of the limitation of {\it smd} to the generalized extreme value distribution appears in the three bias term $\lambda_n$, $\widetilde{\tau}_n$ and $f_{(m)}(x_n) -g_{\bm{\gamma}_m}(x_n)$. $\lambda_n$ corresponds to the asymptotic bias of $\widehat{\bm{\gamma}}_k$ which depends on the second order parameter. The bias $g_{\bm{\gamma}_m}(x_n) -g_{\bm{\gamma}_k}(x_n)$ comes from the asymptotic difference between $m$ and $k$. Though large $k$ yields large variance of $\widehat{\bm{\gamma}}_k$, the precise approximation of $g_{\bm{\gamma}_m}(x_n)$ to $f_{(m)}(x_n)$ (i.e. $\widetilde{\tau}_n$) needs large $m$. As the extreme index tends to zero i.e. $\gamma \to 0$, the convergence rate of the approximation of $g_{\bm{\gamma}_m}(x_n)$ to $f_{(m)}(x_n)$ becomes slow. In the Weibull cases, the rate is slower than any polynomial. Therefore, the convergence rate is sensitive to the first order parameter $\gamma$. Numerical examples in some cases are given later in Table 1 in Section 4.

We next consider {\it smd} estimation without limitation to the generalized extreme value distribution. By using the nonparametric approach, the plug-in type of the kernel density estimator of {\it smd} defined as $\widehat{f}_{(m)} := m \widehat{f}\widehat{F}^{m-1}$ is constructed, where $\widehat{f}$ and $\widehat{F}$ are the kernel density and kernel distribution estimators of $f$ and $F$. We also consider the Block-Maxima-based (BM-based) kernel density estimator $\check{f}_{(m)}$ which is another candidate for the {\it smd} estimator. The BM-based estimator is the kernel density estimator of {\it smd} itself, which is a function of Block-Maxima. From the definition, we see that the nonparametric estimators are consistent for quite wide class of distributions in fixed $m$ cases, where PE is inconsistent. The nonparametric estimators can be good candidate for the {\it smd} at least for relatively small $m$. 

Two nonparametric estimators of {\it smd} aim at the {\it smd} itself directly, and they are free from the misspecification the approximation of $g_{\bm{\gamma}_m}(x_n)$ to $f_{(m)}(x_n)$ (i.e. $\widetilde{\tau}_n$). We theoretically investigate asymptotic properties of the two kernel type of the distribution estimator. Section 2 and 3 gives asymptotic properties of the plug-in kernel density estimator $\widehat{f}_{(m)}$ and the Block-Maxima-based (BM-based) kernel density estimator $\check{f}_{(m)}$. The plug-in kernel density estimator is proved to be theoretically superior to the Block-Maxima-based (BM-based) kernel density estimator. We compare the numerical accuracies of PE $g_{\widehat{\bm{\gamma}}_k}$ and the plug-in estimator $\widehat{f}_{(m)}$ in Section 4. The proofs are given in Appendices.

\section{Nonparametric plug-in type of estimation of {\it smd}}
Applying the plug-in rule, we have the following kernel-type of nonparametric density estimator (NE1):
$$\widehat{f}_{(m)}(x; h_1, h_2) := m \widehat{f}(x; h_1)\widehat{F}^{m-1}(x; h_2)$$
where $\widehat{f}$ and $\widehat{F}$ are the kernel estimator given by
$$\widehat{f}(x; h_1)=\frac{1}{n h_1} \sum_{i=1}^n k\left(\frac{x -X_i}{h_1} \right), ~~~\widehat{F}(x; h_2)=\frac{1}{n} \sum_{i=1}^n K\left(\frac{x -X_i}{h_2} \right),$$
where $k$ is a symmetric and bounded density function whose support is bounded for {\rm (iii)}. $K$ is the cumulative distribution function. $h_j$ are bandwidths that satisfy $n h_1 \to \infty$, $h_j \rightarrow 0$ for $j=1,2$  and $h_j (x^* -x_n) \to 0$ under Assumption 3 for $j=1,2$.
\begin{assumption}
$f$ is twice continuously differentiable at $x_n$.
\end{assumption}
Assumption 8 for the bounded class of distributions requires $\mu$ is integer or $\mu<-3$.
\begin{assumption}
$$\int z^2 k(z) {\rm d}z < \infty, ~~~ \int k^2(z) {\rm d}z < \infty, ~~~ \int z K(z) k(z) {\rm d}z < \infty.$$
\end{assumption}
Using the asymptotic normalities of the estimators, we have the following theorem.
\begin{theorem}
Given Assumption 1, 8 and 9 If $b_{n,1} \to 0$ and $v_{n,1} \to 0$, then $v_{n,1}^{-1/2} \{ (f_{(m)}(x_n) -\widehat{f}_{(m)}(x_n; h_1, h_2)) - b_{n,1}\}$ converges in distribution to the standard normal, where 
\begin{align*}
b_{n,1} :=& \frac{h_1^{2}}{2} \exp(-M_n) M_n^{1+3\gamma} m^{-3\gamma} {\psi_n} \int z^2 k(z) {\rm d}z + \frac{h_2^{2}}{2} M_n m {\xi_n} f(x_n) \int z^2 k(z) {\rm d}z, \\
v_{n,1} :=& \frac{m^2}{nh_1} \exp(-2M_n) f(x_n) \int k^2(z) {\rm d}z + \frac{m^2}{n} M_n^2 f^2(x_n) (M_n^{-1}m - 2h_2 \omega_n \int z K(z) k(z) {\rm d}z), \\
\psi_n:=&\begin{dcases}
\alpha(\alpha+1)(\alpha+2) A^{-3\gamma} ~~~ &\text{for} ~~~ {\rm (i)}\\
\kappa^3 C^3 x_n^{3\kappa-3} ~~~ &\text{for} ~~~ {\rm (ii)}\\
-\mu(\mu+1)(\mu+2) D^{-3\gamma} ~~~ &\text{for} ~~~ {\rm (iii)},
\end{dcases}\\
\xi_n:=&\begin{dcases}
{\alpha (\alpha+1)} x_n^{-2} ~~~ &\text{for} ~~~ {\rm (i)}\\
\kappa^2 C^{2} x_n^{2\kappa-2} ~~~ &\text{for} ~~~ {\rm (ii)}\\
{\mu (\mu +1)} (x^* -x_n)^{-2} ~~~ &\text{for} ~~~ {\rm (iii)}
\end{dcases} \\
\omega_n:=& \begin{dcases}
A^{-1} \alpha x_n^{\alpha-1} ~~~ &\text{for} ~~~ {\rm (i)}\\
\kappa C x_n^{\kappa-1} \exp(C x_n^{\kappa}) ~~~ &\text{for} ~~~ {\rm (ii)}\\
-D^{-1} \mu (x^* -x_n)^{\mu -1} ~~~ &\text{for} ~~~ {\rm (iii)}.
\end{dcases}
\end{align*}
\end{theorem}
\begin{cor}
The asymptotically optimal values of $h_1$ and $h_2$ are given by
$$\left( M_n^{-2-6\gamma} m^{2+6\gamma} \psi_n^{-2} f(x_n) n^{-1} \left(\int z^2 k(z) {\rm d}z\right)^{-2}{\int k^2(z) {\rm d}z}\right)^{1/5}$$
and 
$$\left(2 \xi_n^{-2} \omega_n n^{-1} \left(\int z^2 k(z) {\rm d}z\right)^{-2}{\int z K(z) k(z) {\rm d}z}\right)^{1/3},$$
which are asymptotically identical to those of $\widehat{f}(x_n ; h_1)$ and $\widehat{F}(x_n ; h_2)$. Under the assumption of Theorem 10, $f_{(m)}(x_n) -\widehat{f}_{(m)}(x_n; h_1, h_2)$ with the optimal bandwidths is asymptotically non-degenerate normal with the following asymptotic mean:
$$\nu_1 \exp(-M_n) M_n^{(1+3\gamma)/5} m^{(4-3\gamma)/5} \psi_n^{1/5} f^{2/5}(x_n) n^{-2/5}+ \nu_2 M_n m \xi_n^{-1/3} f(x_n) \omega_n^{2/3} n^{-2/3}$$
where 
\begin{align*}
\nu_1:= \frac{1}{2}\left(\int z^2 k(z) {\rm d}z\right)^{1/5} \left(\int k^2(z) {\rm d}z\right)^{2/5}, \nu_2:= \left(2 \int z^2 k(z) {\rm d}z\right)^{-1/3} \left(\int z K(z) k(z) {\rm d}z\right)^{2/3}.
\end{align*}
\end{cor}
\begin{cor}
The asymptotically optimal bandwidths for $M_n=O(1)$ are of the order
$$\left(\frac{m}{n}\right)^{1/5} \times
\begin{dcases}
m^{\gamma} ~~~ &\text{for} ~~~ {\rm (i), (ii)}\\
\{C^{-1}\ln m \}^{-\theta} ~~~ &\text{for} ~~~ {\rm (ii)}\\
\end{dcases}
$$ 
and
$$\left(\frac{m}{n} \right)^{1/3} \times
\begin{dcases}
m^{\gamma} ~~~ &\text{for} ~~~ {\rm (i), (iii)}\\
\{C^{-1}\ln m \}^{-\theta} ~~~ &\text{for} ~~~ {\rm (ii)}.
\end{dcases}
$$
$\widehat{f}_{(m)}(x_n; h_1, h_2)$ for $M_n=O(1)$ with the optimal bandwidths has the following asymptotic bias of the order:
$$\left(\frac{m}{n}\right)^{2/5}
\begin{dcases}
m^{-\gamma} \\
(\ln m)^{\theta} 
\end{dcases}
+\left(\frac{m}{n}\right)^{2/3} \times
\begin{dcases}
m^{-\gamma} ~~~ &\text{for} ~~~ {\rm (i), (iii)}\\
(\ln m)^{\theta}  ~~~ &\text{for} ~~~ {\rm (ii)}.
\end{dcases}
$$
\end{cor}
The proof of Theorem 10 is found in Appendices. Corollary 11 and 12 are directly obtained from Theorem 10.

\section{Nonparametric block-maxima-based estimation of {\it smd}}
Suppose that the block size is $m$ and ${}^{\exists}\check{n} \in \mathbb{N}$ s.t. $n=\check{n}\times m$. We consider the following BM-based kernel-type of nonparametric density estimator (NE2)
$$\check{f}_{(m)}(x)=\frac{1}{\check{n} h} \sum_{i=1}^{\check{n}} k\left(\frac{x -Y_i}{h} \right),$$
where $Y_j := \max\{X_{m(j-1)+1},X_{m(j-1)+2},\cdots,X_{mj}\}$ $(j=1,\cdots,\check{n})$. The kernel function $k$ is a symmetric and bounded density function whose support is bounded for {\rm (iii)}. $h$ is a bandwidth that satisfies all $h \rightarrow 0$, $\check{n}h \to \infty$ and $h (x^* -x_n) \to 0$ under Assumption 3. The following theorem on the asymptotic normality of the naive nonparametric density estimator holds.
\begin{theorem}
Given Assumption 1, 8 and 9. If $x_n^{\kappa-1} h\to 0$, 
$$b_{n,2} := \frac{mh^2}{2} f(x_n) \phi_n \int z^2 k\left(z \right) {\rm d}z \to 0$$ 
and $v_{n,2} := (\check{n}h)^{-1} f_{(m)}(x_n) \int k^2(z) {\rm d}z \to 0$, then $v_{n,2}^{-1/2} \{ (f_{(m)}(x_n) -\check{f}_{(m)}(x_n; h)) - b_{n,2}\}$ converge in distribution to the standard normal, where 
\begin{align*}
\phi_n := \begin{dcases}
m^2 ~~~ &\text{for} ~~~ {\rm (i), (iii)}\\
(\kappa C)^2 x_n^{2(\kappa-1)} - 2m \kappa C x_n^{\kappa -1} +m^2 ~~~ &\text{for} ~~~ {\rm (ii)}.
\end{dcases}
\end{align*}
\end{theorem}
\begin{cor}
The asymptotically optimal value of $h$ is given by
$$\left( n^{-1} \{f(x_n)\}^{-1} \phi_n^{-2} \left(\int z^2 k(z) {\rm d}z\right)^{-2}{\int k^2(z) {\rm d}z}\right)^{1/5}.$$ 
Under the assumption of Theorem 13, $f_{(m)}(x_n) -\check{f}_{(m)}(x_n; h)$ with the optimal bandwidths is asymptotically non-degenerate normal with the following asymptotic mean:
$$\nu_1 n^{-2/5} m \{f(x_n)\}^{3/5} \phi_n^{1/5}.$$
\end{cor}
\begin{cor}
The asymptotically optimal bandwidth for $M_n=O(1)$ is of the order
$$n^{-1/5} m^{-3/5} \times
\begin{dcases}
m^{\gamma/5} ~~~ &\text{for} ~~~ {\rm (i), (iii)}\\
(\ln m)^{-1/5}  ~~~ &\text{for} ~~~ {\rm (ii)}.
\end{dcases}
$$ 
Under the assumption of Theorem 13, $\check{f}_{(m)}(x_n; h)$ has the asymptotic bias of the order:
$$\left(\frac{m^2}{n}\right)^{2/5}
\begin{dcases}
m^{-(3/5)\gamma} ~~~ &\text{for} ~~~ {\rm (i), (iii)}\\
(\ln m)^{-(3/5)\gamma} ~~~ &\text{for} ~~~ {\rm (ii)}.
\end{dcases}
$$
\end{cor}
Theorem 13 is proved in the Appendices. Corollaries 14 and 15 are the direct consequences. 

\section{Comparative study}
Performances of the three estimators PE, NE1 and NE2 are compared in this section. The first is theoretically and the second is numerically. $M_n \equiv K_n \equiv \delta$ is supposed throughout this section. NE1 and NE2 with their respective optimal bandwidths converge on the order of $m^{-2\gamma}(m/n)^{4/5}$ and $m^{-(6/5)\gamma}(m^2/n)^{4/5}$. The following Theorem 16 is immediately seen from Theorems 10 and 13.

\begin{theorem}
The convergence rate of NE1 $\widehat{f}_{(m)}(x_n; h_1, h_2)$ with the optimal bandwidths is faster than NE2 $\check{f}_{(m)}(x_n)$ with the optimal bandwidth for $\gamma>-1$.
\end{theorem}

PE converges with the rate $m^{-2\gamma} \times (N^{-1} m^{2-4\beta} + m^{-2\gamma\beta} + m^{-2} + N^{-1})$ for the Hall class. For the Weibull class, PE does not satisfy  both Assumption 2 and 4 with any block size $k$ and is inconsistent. For the bounded class, PE has the convergence rate $m^{-2\gamma} \times (N^{-1} m^{2+4\sigma} + m^{-2\gamma\sigma} + m^{-2} + N^{-1})$. Assumption 5 requires $\gamma>-3/2$,  $\gamma>-1/2$ and $\gamma>-1/6$ for $m=n^{1/4}$, $m=n^{1/2}$ and $m=n^{3/4}$ respectively. The parameters of the distributions and the convergence rates of the estimators without terms slower than any polynomial are summarized in Table 1. The hyphens mean the distribution breaks the assumption of the theorems. The settings in the density estimation and distribution estimation surveyed in Moriyama [4] are same. The difference is all the target values $f_{(m)}(x_n)$ tend to zero with the exception of the bounded class with $\mu=-1$ while all $F^m(x_n) \to \exp(-\delta)>0$ in distribution estimation. The convergence rates divided by $f_{(m)}(x_n)$ are given in Table 2, where the hyphens are also used to mean that the estimators are inconsistent. All the convergence rates of NE2 are slower than or same as those of NE1 in the Table 2 as Theorem 16 stated. While for $m=n^{1/4}$ all the rates of NE1 are faster than PE, the differences become small as $m$ gets large, which coincides with the results in distribution estimation given in Moriyama [4]. Both the rates of PE and NE1 become slow as $\gamma$ gets small. 

By simulating the mean integrated squared error ($MISE$) of PE $g_{\widehat{\bm{\gamma}}_k}$:
$$L_m \int_{Q_ m(0.1)}^{Q_ m(0.9)} \left(g_{\widehat{\bm{\gamma}}_k}(x) -f_{(m)}(x)\right)^2 {\rm d}x.$$ 
and that of $\widehat{f}_{(m)}$, we studied the numerical accuracy in small-sample cases, where $L_m := Q_ m(0.9) - Q_ m(0.1)$ and $Q_ m(q)$ denotes the $q$th quantile of the {\it smd}. Tables 3--4 show the mean values and their standard deviation (sd) of the obtained $MISE$ values by simulating 1000 times. The underlying distributions $F$ were Pareto distributions, $T$-distributions, Burr distributions, Fr\'{e}chet distributions, Weibull distributions and reversed Burr distributions. The forecast periods were $m =n^{1/4}$, $m=n^{1/2}$, and $m=n^{3/4}$. All kernel functions were the Epanechnikov for the bounded class and the Gaussian for the other classes. Both of the bandwidths $h_1$ and $h_2$ were estimated by the cross-validation approach or the plug-in approach in the Tables `CV' and `PI' respectively. CV means the unbiased cross-validation estimator (see e.g. Silverman [8]) and Bowman et al [9]'s cross-validation estimator. As the plug-in approaches, the Sheather and Jones [10]'s method and the Altman and Leger [11]'s method were employed. The sample sizes were $(n=)2^{8}$ in Tables 3 and $(n=)2^{12}$ in Tables 4 . 

Comparing CV with PI, we see CV greatly outperforms for $1 \le \gamma \le 2$ while PI is much better for $\gamma \ge 4$. In the other cases, they are numerically comparable. Comparing NE1 with PE, we find that the MISE values of PE were smaller in most of the cases with $m=n^{1/4}$. For $\gamma \fallingdotseq 0$ those of NE1 were smaller as $m$ get large. For the light-tailed distributions with $\gamma<0$, the MISE values of NE1 were comparable with those of PE, though so-called boundary bias in NE1 for bounded distributions with $\gamma \ge -1$ was observed. Due to the boundary bias, the results especially for distributions with large extreme indices (e.g. $\gamma \ge -1/3$) are different from those of distribution estimation shown in Moriyama [4]. Comparing the cases $(n=)2^{12}$ with $(n=)2^{8}$, it is seen that PE substantially reduces the MISE values, which are smaller than those of NE1.

\begin{landscape}
\begin{table}
\caption{The polynomial convergence rates of MSE of the estimators}{\fontsize{6pt}{6pt}\selectfont
$$\begin{tabu}[c]{|ccc||ccc|ccc|ccc|}
  \hline
 \multicolumn{3}{|c||}{{\rm Pareto}} & \multicolumn{3}{c|}{m=n^{1/4}} & \multicolumn{3}{c|}{m=n^{1/2}} & \multicolumn{3}{c|}{m=n^{3/4}} \\
\ell & \alpha & \beta & {\rm PE} & {\rm NE1} & {\rm NE2} & {\rm PE} & {\rm NE1} & {\rm NE2} & {\rm PE} & {\rm NE1} & {\rm NE2} \\ \hline
 1/2 & 1/2 & 1 & -1/2 & -8/5 & -1 & -1 & -12/5 & -6/5 & -3/2 & -16/5 & -7/5 \\ 
 1 & 1 & 1 & -1/2& -11/10 & -7/10 & -1 & -7/5 & -3/5 & -3/2 & -17/10 & -1/2 \\ 
 3 & 3 & 1 & -1/3 & -23/30 & -1/2 & -2/3 & -11/15 & -1/5 & -3/4 & -7/10 & \text{--} \\ 
 10 & 10 & 1 & -1/10 & -13/20 & -43/100 & -1/5 & -1/2 & -3/50 & -3/10 & -7/20 & \text{--} \\ 
 
\hline\hline
 \multicolumn{3}{|c||}{{\rm T} dist.} & \multicolumn{3}{c|}{m=n^{1/4}} & \multicolumn{3}{c|}{m=n^{1/2}} & \multicolumn{3}{c|}{m=n^{3/4}} \\
\ell & \alpha & \beta & {\rm PE} & {\rm NE1} & {\rm NE2} & {\rm PE} & {\rm NE1} & {\rm NE2} & {\rm PE} & {\rm NE1} & {\rm NE2} \\ \hline
 1/2 & 1/2 & 2 & -3/2 & -8/5 & -1 & -5/2 & -12/5 & -6/5 & -13/4 & -16/5 & -7/5 \\ 
 1 & 1 & 2 & -1/2& -11/10 & -7/10 & -3/2 & -7/5 & -3/5 & -7/4 & -17/10 & -1/2 \\ 
 3 & 3 & 2 & -1/3 & -23/30 & -1/2 & -2/3 & -11/15 & -1/5 & -3/4 & -7/10 & \text{--} \\ 
 10 & 10 & 2 & -1/10 & -13/20 & -43/100 & -1/5 & -1/2 & -3/50 & -3/10 & -7/20 & \text{--} \\ 
 \hline\hline
  \multicolumn{3}{|c||}{{\rm Burr}} & \multicolumn{3}{c|}{m=n^{1/4}} & \multicolumn{3}{c|}{m=n^{1/2}} & \multicolumn{3}{c|}{m=n^{3/4}} \\
c,\ell & \alpha & \beta & {\rm PE} & {\rm NE1} & {\rm NE2} & {\rm PE} & {\rm NE1} & {\rm NE2} & {\rm PE} & {\rm NE1} & {\rm NE2} \\ \hline
 1/2,1/2 & 1/4 & 1/2 & -5/2 & -13/5 & -8/5 & -9/2 & -22/5 & -12/5 & -25/4 & -31/5 & -16/5 \\ 
 1,1/2 & 1/2 & 1 & -1/2 & -8/5 & -1 & -1 & -12/5 & -6/5 & -3/2 & -16/5 & -7/5  \\ 
 3,1/2 & 3/2 & 3 & -2/3 & -14/15 & -3/5 & -7/6 & -16/15 & -2/5 & -5/4 & -6/5 & -1/5 \\ 
 1/2,1 & 1/2 & 1/2 & -3/2 & -8/5 & -1 & -3/2 & -12/5 & -6/5 & -13/4 & -16/5 & -7/5 \\ 
 1,1 & 1 & 1 & -1/2 & -11/10 & -7/10 & -1 & -7/5 & -3/5 & -3/2 & -17/10 & -1/2 \\  
 3,1 & 3 & 3 & -1/3 & -23/30 & -1/2 & -2/3 & -11/15 & -1/5 & -3/4 & -7/10 & \text{--} \\ 
 1/2,3 & 3/2 & 1/2 & -1/2 & -14/15 & -3/5 & -1 & -16/15 & -2/5 & -5/4 & -6/5 & -1/5 \\ 
 1,3 & 3 & 1 & -1/3 & -23/30 & -1/2 & -2/3 & -11/15 & -1/5 & -3/4 & -7/10 & \text{--} \\ 
 3,3 & 9 & 3 & -1/9 & -59/90 & -13/30 & -2/9 & -23/45 & -1/15 & -1/3 & -11/30 & \text{--} \\ 
 
  \hline\hline
 \multicolumn{3}{|c||}{{\rm Fr\acute{e}chet}} & \multicolumn{3}{c|}{m=n^{1/4}} & \multicolumn{3}{c|}{m=n^{1/2}} & \multicolumn{3}{c|}{m=n^{3/4}} \\
\gamma & \alpha & \beta & {\rm PE} & {\rm NE1} & {\rm NE2} & {\rm PE} & {\rm NE1} & {\rm NE2} & {\rm PE} & {\rm NE1} & {\rm NE2} \\ \hline
 5 & 1/5 & 1/5 & \text{--} & -31/10 & -19/10 & \text{--} & -27/5 & -3 & \text{--} & -77/10 & -11/10 \\ 
 2 & 1/2 & 1/2 & -3/2 & -8/5 & -1 & -3/2 & -12/5 & -6/5 & -13/4 & -16/5 & -7/5 \\  
 1 & 1 & 1 & -1/2 & -11/10 & -7/10 & -1 & -7/5 & -3/5 & -3/2 & -17/10 & -1/2 \\ 
 1/2 & 2 & 1 & -1/2 & -17/20 & -11/20 & -1 & -9/10 & -3/10 & -1 & -19/20 & -1/20 \\ 
 1/4 & 4 & 1 & -1/4 & -29/40 & -19/40 & -1/2 & -13/20 & -3/20 & -5/8 & -23/40 & \text{--} \\
 
 \hline\hline
 \multicolumn{3}{|c||}{Weibull} & \multicolumn{3}{c|}{m=n^{1/4}} & \multicolumn{3}{c|}{m=n^{1/2}} & \multicolumn{3}{c|}{m=n^{3/4}} \\
\kappa& \gamma & \rho & {\rm PE} & {\rm NE1} & {\rm NE2} & {\rm PE} & {\rm NE1} & {\rm NE2} & {\rm PE} & {\rm NE1} & {\rm NE2} \\ \hline
1/2 & 0 & 0 & \text{--} & -3/5 & -2/5 & \text{--} & -2/5 & \text{--} & \text{--} & -1/5 & \text{--} \\
1 & 0 & 0 & \text{--} & -3/5 & -2/5 & \text{--} & -2/5 & \text{--} & \text{--} & -1/5 & \text{--} \\
3 & 0 & 0 & \text{--} & -3/5 & -2/5 & \text{--} & -2/5 & \text{--} & \text{--} & -1/5 & \text{--} \\
10 & 0 & 0 & \text{--} & -3/5 & -2/5 & \text{--} & -2/5 & \text{--} & \text{--} & -1/5 & \text{--} \\
 \hline\hline
  \multicolumn{3}{|c||}{{\rm rev. Burr}} & \multicolumn{3}{c|}{m=n^{1/4}} & \multicolumn{3}{c|}{m=n^{1/2}} & \multicolumn{3}{c|}{m=n^{3/4}} \\
c,\ell & \mu & \sigma & {\rm PE} & {\rm NE1} & {\rm NE2} & {\rm PE} & {\rm NE1} & {\rm NE2} & {\rm PE} & {\rm NE1} & {\rm NE2} \\ \hline
 -1/2,-1/3 & -6 & -2 & -1/12 & -31/60 & -7/20 & -1/6 & -7/30 & \text{--} & \text{--} & \text{--} & \text{--} \\
 -1,-1/3 & -3 & -1 & \text{--} & -13/30 & -3/10 & \text{--} & -1/15 & \text{--} & \text{--} & \text{--} & \text{--}\\
 -3,-1/3 & -1 & -1/3 & \text{--} & -1/10 & -1/10 & \text{--} & \text{--} & \text{--} & \text{--} & \text{--} & \text{--} \\
 -1/2,-1 & -2 & -2 & \text{--} & -7/20 & -1/4 & \text{--} & \text{--} & \text{--} & \text{--} & \text{--} & \text{--} \\
 -1,-1 & -1 & -1 & \text{--} & -1/10 & -1/10 & \text{--} & \text{--} & \text{--} & \text{--} & \text{--} & \text{--} \\
 -3,1 & -1/3 & -1/3 & \text{--} & \text{--} & \text{--} & \text{--} & \text{--} & \text{--} & \text{--} & \text{--} & \text{--} \\
 -1/2,-2 & -1 & -2 & \text{--} & -1/10 & -1/10 & \text{--} & \text{--} & \text{--} & \text{--} & \text{--} & \text{--}\\
 -1,-2 & -1/2 & -1 & \text{--} & \text{--} & \text{--} & \text{--} & \text{--} & \text{--} & \text{--} & \text{--} & \text{--} \\
 -3,-2 & -1/6 & -1/3 & \text{--} & \text{--} & \text{--} & \text{--} & \text{--} & \text{--} & \text{--} & \text{--} & \text{--} \\ 
 \hline
\end{tabu}$$
}
\end{table}
\end{landscape}

\begin{landscape}
\begin{table}
\caption{The normalized polynomial convergence rates of MSE of the estimators}{\fontsize{6pt}{6pt}\selectfont
$$\begin{tabu}[c]{|ccc||ccc|ccc|ccc|}
  \hline
 \multicolumn{3}{|c||}{{\rm Pareto}} & \multicolumn{3}{c|}{m=n^{1/4}} & \multicolumn{3}{c|}{m=n^{1/2}} & \multicolumn{3}{c|}{m=n^{3/4}} \\
\ell & \alpha & \beta & {\rm PE} & {\rm NE1} & {\rm NE2} & {\rm PE} & {\rm NE1} & {\rm NE2} & {\rm PE} & {\rm NE1} & {\rm NE2} \\ \hline
 1/2 & 1/2 & 1 & \text{--} & -11/10 & -1/2 & \text{--} & -7/5 & -1/5 & \text{--} & -17/10 & \text{--} \\ 
 1 & 1 & 1 & -1/4 & -17/20 & -9/20 & -1/2 & -9/10 & -1/10 & -3/4 & -19/20 & \text{--} \\ 
 3 & 3 & 1 & -1/4 & -41/60 & -5/12 & -1/2 & -17/30 & -1/30 & -1/2 & -9/20 & \text{--} \\ 
 10 & 10 & 1 & -3/40 & -5/8 & -81/200 & -3/20 & -9/20 & -1/100 & -9/40 & -11/40 & \text{--} \\ 
 
\hline\hline
 \multicolumn{3}{|c||}{{\rm T} dist.} & \multicolumn{3}{c|}{m=n^{1/4}} & \multicolumn{3}{c|}{m=n^{1/2}} & \multicolumn{3}{c|}{m=n^{3/4}} \\
\ell & \alpha & \beta & {\rm PE} & {\rm NE1} & {\rm NE2} & {\rm PE} & {\rm NE1} & {\rm NE2} & {\rm PE} & {\rm NE1} & {\rm NE2} \\ \hline
 1/2 & 1/2 & 2 & -1 & -11/10 & -1/2 & -3/2 & -7/5 & -1/5 & -7/4 & -17/10 & \text{--} \\ 
 1 & 1 & 1 & -1/4 & -17/20 & -9/20 & -1 & -9/10 & -1/10 & -1 & -19/20 & \text{--} \\ 
 3 & 3 & 1 & -1/4 & -41/60 & -5/12 & -1/2 & -17/30 & -1/30 & -1/2 & -9/20 & \text{--} \\ 
 10 & 10 & 1 & -3/40 & -5/8 & -81/200 & -3/20 & -9/20 & -1/100 & -9/40 & -11/40 & \text{--} \\ 
 
 \hline\hline
  \multicolumn{3}{|c||}{{\rm Burr}} & \multicolumn{3}{c|}{m=n^{1/4}} & \multicolumn{3}{c|}{m=n^{1/2}} & \multicolumn{3}{c|}{m=n^{3/4}} \\
c,\ell & \alpha & \beta & {\rm PE} & {\rm NE1} & {\rm NE2} & {\rm PE} & {\rm NE1} & {\rm NE2} & {\rm PE} & {\rm NE1} & {\rm NE2} \\ \hline
 1/2,1/2 & 1/4 & 1/2 & -3/2 & -8/5 & -3/5 & -5/2 & -12/5 & -2/5 & -13/4 & -16/5 & -1/5 \\ 
 1,1/2 & 1/2 & 1 & \text{--} & -11/10 & -1/2 & \text{--} & -7/5 & -1/5 & \text{--} & -17/10 & \text{--} \\ 
 3,1/2 & 3/2 & 3 & -7/24 & -67/120 & -9/40 & -3/28 & -19/60 & \text{--} & -1/8 & -3/40 & \text{--} \\ 
 1/2,1 & 1/2 & 1/2 & -1 & -11/10 & -1/2 & -1/2 & -7/5 & -1/5 & -7/4 & -17/10 & \text{--} \\ 
 1,1 & 1 & 1 & -1/4 & -17/20 & -9/20 & -1 & -9/10 & -1/10 & -1 & -19/20 & \text{--} \\ 
 3,1 & 3 & 3 & -1/4 & -41/60 & -5/12 & -1/2 & -17/30 & -1/30 & -1/2 & -9/20 & \text{--} \\ 
 1/2,3 & 3/2 & 1/2 & -1/8 & -67/120 & -9/40 & -1/4 & -19/60 & \text{--} & -1/8 & -3/40 & \text{--} \\ 
 1,3 & 3 & 1 & -1/4 & -41/60 & -5/12 & -1/2 & -17/30 & -1/30 & -1/2 & -9/20 & \text{--} \\ 
 3,3 & 9 & 3 & -1/12 & -113/180 & -73/180 & -1/6 & -41/90 & -1/90 & -1/4 & -17/60 & \text{--} \\ 
 
  \hline\hline
 \multicolumn{3}{|c||}{{\rm Fr\acute{e}chet}} & \multicolumn{3}{c|}{m=n^{1/4}} & \multicolumn{3}{c|}{m=n^{1/2}} & \multicolumn{3}{c|}{m=n^{3/4}} \\
\gamma & \alpha & \beta & {\rm PE} & {\rm NE1} & {\rm NE2} & {\rm PE} & {\rm NE1} & {\rm NE2} & {\rm PE} & {\rm NE1} & {\rm NE2} \\ \hline
 5 & 1/5 & 1/5 & \text{--} & -74/40 & -26/40 & \text{--} & -29/10 & -1/2 & \text{--} & -79/20 & \text{--} \\ 
 2 & 1/2 & 1/2 & -1 & -11/10 & -1/2 & -1/2 & -7/5 & -1/5 & -7/4 & -17/10 & \text{--} \\ 
 1 & 1 & 1 & -1/4 & -17/20 & -9/20 & -1 & -9/10 & -1/10 & -1 & -19/20 & \text{--} \\ 
 1/2 & 2 & 1 & -3/8 & -29/40 & -17/40 & -3/4 & -13/20 & -1/20 & -5/8 & -21/40 & \text{--} \\ 
 1/4 & 4 & 1 & -3/16 & -53/80 & -33/80 & -3/8 & -21/40 & -1/40 & -7/16 & -31/80 & \text{--} \\
 
 \hline\hline
 \multicolumn{3}{|c||}{Weibull} & \multicolumn{3}{c|}{m=n^{1/4}} & \multicolumn{3}{c|}{m=n^{1/2}} & \multicolumn{3}{c|}{m=n^{3/4}} \\
\kappa& \gamma & \rho & {\rm PE} & {\rm NE1} & {\rm NE2} & {\rm PE} & {\rm NE1} & {\rm NE2} & {\rm PE} & {\rm NE1} & {\rm NE2} \\ \hline
1/2 & 0 & 0 & \text{--} & -3/5 & -2/5 & \text{--} & -2/5 & \text{--} & \text{--} & -1/5 & \text{--} \\
1 & 0 & 0 & \text{--} & -3/5 & -2/5 & \text{--} & -2/5 & \text{--} & \text{--} & -1/5 & \text{--} \\
3 & 0 & 0 & \text{--} & -3/5 & -2/5 & \text{--} & -2/5 & \text{--} & \text{--} & -1/5 & \text{--} \\
10 & 0 & 0 & \text{--} & -3/5 & -2/5 & \text{--} & -2/5 & \text{--} & \text{--} & -1/5 & \text{--} \\
 \hline\hline
  \multicolumn{3}{|c||}{{\rm rev. Burr}} & \multicolumn{3}{c|}{m=n^{1/4}} & \multicolumn{3}{c|}{m=n^{1/2}} & \multicolumn{3}{c|}{m=n^{3/4}} \\
c,\ell & \mu & \sigma & {\rm PE} & {\rm NE1} & {\rm NE2} & {\rm PE} & {\rm NE1} & {\rm NE2} & {\rm PE} & {\rm NE1} & {\rm NE2} \\ \hline
 -1/2,-1/3 & -6 & -2 & -1/8 & -67/120 & -47/120 & -1/4 & -19/60 & \text{--} & \text{--} & \text{--} & \text{--} \\
 -1,-1/3 & -3 & -1 & \text{--} & -31/60 & -23/60 & \text{--} & -7/30 & \text{--} & \text{--} & \text{--} & \text{--}\\
 -3,-1/3 & -1 & -1/3 & \text{--} & -7/20 & -7/20 & \text{--} & \text{--} & \text{--} & \text{--} & \text{--} & \text{--} \\
 -1/2,-1 & -2 & -2 & \text{--} & -19/40 & -3/8 & \text{--} & \text{--} & \text{--} & \text{--} & \text{--} & \text{--} \\
 -1,-1 & -1 & -1 & \text{--} & -7/20 & -7/20 & \text{--} & \text{--} & \text{--} & \text{--} & \text{--} & \text{--} \\
 -3,1 & -1/3 & -1/3 & \text{--} & \text{--} & \text{--} & \text{--} & \text{--} & \text{--} & \text{--} & \text{--} & \text{--} \\
 -1/2,-2 & -1 & -2 & \text{--} & -7/20 & -7/20 & \text{--} & \text{--} & \text{--} & \text{--} & \text{--} & \text{--}\\
 -1,-2 & -1/2 & -1 & \text{--} & \text{--} & \text{--} & \text{--} & \text{--} & \text{--} & \text{--} & \text{--} & \text{--} \\
 -3,-2 & -1/6 & -1/3 & \text{--} & \text{--} & \text{--} & \text{--} & \text{--} & \text{--} & \text{--} & \text{--} & \text{--} \\ 
 \hline
\end{tabu}$$
}
\end{table}
\end{landscape}

\begin{landscape}
\begin{table}
\caption{Scaled MISE values and sd values where $n=2^8$}{\fontsize{6pt}{6pt}\selectfont
$$\begin{tabu}[c]{|c||cccccc|cccccc|cccccc|}
   \hline
    {\rm Pareto} & \multicolumn{6}{c|}{m=n^{1/4}} & \multicolumn{6}{c|}{m=n^{1/2}} & \multicolumn{6}{c|}{m=n^{3/4}} \\
 \ell & {\rm PE} & {\rm sd} & {\rm CV} & {\rm sd}  & {\rm PI} & {\rm sd}  & {\rm PE} & {\rm sd} & {\rm CV} & {\rm sd} & {\rm PI} & {\rm sd} & {\rm PE} & {\rm sd} & {\rm CV} & {\rm sd} & {\rm PI} & {\rm sd} \\ \hline
1/2 & 0.242 & 0.443 & 8.065 & 1.721 & 7.059 & 11.99 & 9.222 & 108.9 & 9.339 & 3.345 & 487.7 & 6491 & 51.01 & 383.8 & 108.6 & 296.2 & 22770 &  311400 \\ 
  1 & 0.037 & 0.041 & 0.353 & 0.393 & 0.771 & 0.422 & 0.510 & 5.696 & 1.474 & 1.156 & 7.401 & 12.05 & 33.36 & 696.2 & 20.49 & 27.17 & 122.8 & 500.6 \\ 
  3 & 0.019 & 0.021 & 0.133 & 0.057 & 0.088 & 0.036 & 0.116 & 0.162 & 0.799 & 0.344 & 0.514 & 0.178 & 54.79 & 1202 & 5.143 & 2.904 & 3.231 & 1.437 \\ 
  10 & 0.016 & 0.017 & 0.102 & 0.046 & 0.054 & 0.025 & 0.104 & 0.176 & 0.447 & 0.202 & 0.240 & 0.105 & 7.608 & 71.76 & 2.073 & 1.007 & 1.093 & 0.484 \\ 
  
\hline\hline
 {\rm T} dist. & \multicolumn{6}{c|}{m=n^{1/4}} & \multicolumn{6}{c|}{m=n^{1/2}} & \multicolumn{6}{c|}{m=n^{3/4}} \\
\ell & {\rm PE} & {\rm sd} & {\rm CV} & {\rm sd}  & {\rm PI} & {\rm sd}  & {\rm PE} & {\rm sd} & {\rm CV} & {\rm sd} & {\rm PI} & {\rm sd} & {\rm PE} & {\rm sd} & {\rm CV} & {\rm sd} & {\rm PI} & {\rm sd} \\ \hline
1/2 & 0.996 & 0.986 & 7.695 & 0.921 & 4.380 & 2.679 & 6.928 & 53.80 & 9.947 & 0.953 & 146.1 & 3708 & 41.92 & 205.4 & 31.32 & 61.64 & 28650 & 856500 \\ 
  1 & 0.068 & 0.067 & 0.279 & 0.323 & 0.314 & 0.267 & 0.367 & 0.744 & 0.909 & 0.446 & 2.751 & 7.268 & 42.09 & 695.3 & 8.641 & 9.806 & 45.99 & 296.1 \\ 
  3 & 0.016 & 0.014 & 0.021 & 0.020 & 0.018 & 0.014 & 0.101 & 0.213 & 0.095 & 0.075  & 0.088 & 0.065 & 12.73 & 215.9 & 0.507 & 0.336 & 0.496 & 0.322 \\ 
  10 & 0.011 & 0.012 & 0.015 & 0.015 & 0.012 & 0.010 & 0.088 & 0.137 & 0.049 & 0.044 & 0.042 & 0.036 & 7.153 & 72.64 & 0.193 & 0.177 & 0.162 & 0.144 \\ 
  
 \hline\hline
 {\rm Burr} & \multicolumn{6}{c|}{m=n^{1/4}} & \multicolumn{6}{c|}{m=n^{1/2}} & \multicolumn{6}{c|}{m=n^{3/4}} \\
 c,\ell & {\rm PE} & {\rm sd} & {\rm CV} & {\rm sd}  & {\rm PI} & {\rm sd}  & {\rm PE} & {\rm sd} & {\rm CV} & {\rm sd} & {\rm PI} & {\rm sd} & {\rm PE} & {\rm sd} & {\rm CV} & {\rm sd} & {\rm PI} & {\rm sd} \\ \hline
1/2,1/2 & 102.4 & 156.6 & 714.1 & 1.879 & 361.2 & 139.4 & 548.3 & 1212 & 942.8 & 14.97 & 942.3 & 0.909 & 1306 & 6052 & 9436 & 75610 & 1037 & 0.000 \\ 
  1,1/2 & 0.217 & 0.416 & 8.115 & 1.693 & 7.056 & 11.85 & 5.507 & 22.83 & 9.281 & 3.120 & 385.1 & 5181 & 43.43 & 360.3 & 104.3 & 258.9 & 84390 & 2319000 \\ 
  3,1/2 & 0.023 & 0.029 & 0.097 & 0.068 & 0.113 & 0.093 & 0.178 & 0.341 & 0.754 & 0.386 & 0.866 & 0.485 & 12.23 & 107.1 & 7.402 & 5.880 & 8.760 & 9.997 \\ 
  1/2,1 & 0.411 & 0.674 & 11.52 & 1.685 & 10.20 & 27.86 & 7.554 & 46.29 & 10.36 & 4.495 & 1590 & 41270 & 44.30 & 288.4 & 117.4 & 394.4 & 11970 & 212600 \\ 
  1,1 & 0.041 & 0.046 & 0.309 & 0.350 & 0.849 & 0.470 & 0.391 & 1.235 & 1.456 & 1.274 & 8.263 & 12.39 & 17.72 & 160.2  & 20.57 & 26.62 & 159.3 & 738.4 \\ 
  3,1 & 0.014 & 0.015 & 0.034 & 0.025 & 0.031 & 0.019 & 0.112 & 0.190 & 0.187 & 0.121 & 0.176 & 0.090 & 7.005 & 39.81 & 1.117 & 0.676 & 1.072 & 0.509 \\ 
  1/2,3 & 0.072 & 0.072 & 0.242 & 0.209 & 1.355 & 0.698 & 0.291 & 1.140 & 1.266 & 0.976 & 15.32 & 23.85 & 53.22 & 1205 & 13.95 & 14.99 & 178.7 & 748.7 \\ 
  1,3 & 0.019 & 0.021 & 0.132 & 0.055 & 0.085 & 0.035 & 0.123 & 0.202 & 0.804 & 0.343 & 0.504 & 0.173 & 81.56 & 1379 & 5.059 & 2.907 & 3.255 & 1.361 \\ 
  3,3 & 0.013 & 0.015 & 0.017 & 0.018 & 0.015 & 0.012 & 0.091 & 0.167 & 0.065 & 0.065 & 0.056 & 0.043 & 18.88 & 437.1 & 0.263 & 0.235 & 0.232 & 0.182 \\ 
  
   \hline\hline
  {\rm Fr\acute{e}chet} & \multicolumn{6}{c|}{m=n^{1/4}} & \multicolumn{6}{c|}{m=n^{1/2}} & \multicolumn{6}{c|}{m=n^{3/4}} \\
 \gamma & {\rm PE} & {\rm sd} & {\rm CV} & {\rm sd}  & {\rm PI} & {\rm sd}  & {\rm PE} & {\rm sd} & {\rm CV} & {\rm sd} & {\rm PI} & {\rm sd} & {\rm PE} & {\rm sd} & {\rm CV} & {\rm sd} & {\rm PI} & {\rm sd} \\ \hline
5 & 4449 & 4841 & 13700 & 1.368 & 5955 & 7608 & 11340 & 20670 & 13710 & 13690 & 29.60 & 108.4 & 13990 & 16200 & 74840 & 722000 & 13700 & 0.000 \\ 
  2 & 0.303 & 0.741 & 9.454 & 1.692 & 8.314 & 12.75 & 6.025 & 38.12 & 9.718 & 3.716 & 235.7 & 2446 & 61.86 & 412.0 & 114.2 & 316.4 & 74590 & 1840000 \\ 
  1 & 0.041 & 0.054 & 0.329 & 0.374 & 0.829 & 0.561 & 0.381 & 1.065 & 1.505 & 1.232 & 9.181 & 21.77 & 123.3 & 3281 & 20.28 & 27.15 & 246.2 & 2078 \\ 
  1/2 & 0.019 & 0.022 & 0.082 & 0.043 & 0.076 & 0.042 & 0.139 & 0.208 & 0.588 & 0.281 & 0.519 & 0.184 & 20.72 & 343.2 & 4.717 & 3.338 &  3.991 & 1.861 \\ 
  1/4 & 0.015 & 0.018 & 0.040 & 0.027 & 0.035 & 0.021 & 0.102 & 0.158 & 0.197 & 0.125 & 0.176 & 0.096 & 47.66 & 664.3 & 1.114 & 0.692 & 0.972 & 0.458 \\ 
      
  \hline\hline
    {\rm Weibull} & \multicolumn{6}{c|}{m=n^{1/4}} & \multicolumn{6}{c|}{m=n^{1/2}} & \multicolumn{6}{c|}{m=n^{3/4}} \\
 \kappa & {\rm PE} & {\rm sd} & {\rm CV} & {\rm sd}  & {\rm PI} & {\rm sd}  & {\rm PE} & {\rm sd} & {\rm CV} & {\rm sd} & {\rm PI} & {\rm sd} & {\rm PE} & {\rm sd} & {\rm CV} & {\rm sd} & {\rm PI} & {\rm sd} \\ \hline
1/2 & 0.044 & 0.038 & 0.159 & 0.067 & 0.299 & 0.085 & 0.108 & 0.169 & 0.978 & 0.438 & 1.828 & 0.702 & 54.80 & 1520 & 5.554 & 2.943 & 10.18 & 7.603 \\ 
  1 & 0.017 & 0.020 & 0.085 & 0.045 & 0.044 & 0.024 & 0.092 & 0.148 & 0.333 & 0.178 & 0.167 & 0.087 & 39.12 & 621.0 & 1.285 & 0.712 & 0.669 & 0.343 \\ 
  3 & 0.011 & 0.013 & 0.012 & 0.013 & 0.010 & 0.009 & 0.090 & 0.171 & 0.035 & 0.032 & 0.031 & 0.027 & 18.64 & 242.2 & 0.117 & 0.099 & 0.112 & 0.101 \\ 
  10 & 0.011 & 0.012 & 0.012 & 0.011 & 0.009 & 0.008 & 0.088 & 0.127 & 0.034 & 0.030 & 0.027 & 0.023 & 8.291 & 127.5 & 0.126 & 0.091 & 0.098 & 0.073 \\ 
  
   \hline\hline
 {\rm rev. Burr} & \multicolumn{6}{c|}{m=n^{1/4}} & \multicolumn{6}{c|}{m=n^{1/2}} & \multicolumn{6}{c|}{m=n^{3/4}} \\
 c,\ell & {\rm PE} & {\rm sd} & {\rm CV} & {\rm sd}  & {\rm PI} & {\rm sd}  & {\rm PE} & {\rm sd} & {\rm CV} & {\rm sd} & {\rm PI} & {\rm sd} & {\rm PE} & {\rm sd} & {\rm CV} & {\rm sd} & {\rm PI} & {\rm sd} \\ \hline
 -1/2,-1/3 & 24500 & 0.000 & 24500 & 0.000 & 42700 & 0.000 & 54800 & 0.000 & 42700 & 0.000 & 33300 & 0.000 & 54800 & 0.000 & 54800 & 0.000 & 59600 & 0.000 \\ 
  -1,-1/3 & 56.10 & 0.203 & 56.60 & 0.000 & 30.70 & 0.000 & 26.50 & 0.000 & 30.70 & 0.000 & 37.40 & 0.000 & 26.50 & 0.000 & 26.50 & 0.000 & 25.40 & 0.000 \\ 
 -3,-1/3 & 4.675 & 0.129 & 5.270 & 0.000 & 1.890 & 0.000 & 1.300 & 0.000 & 1.890 & 0.000 & 2.718 & 0.004 & 1.300 & 0.000 & 1.300 & 0.000 & 1.040 & 0.000 \\ 
 -1/2,-1 & 1.524 & 0.184 & 2.400 & 0.192 & 3.027 & 0.022 & 3.260 & 0.000 & 3.043 & 0.014 & 1.921 & 0.471 & 3.260 & 0.000 & 3.260 & 0.000 & 3.340 & 0.000 \\ 
  -1,-1 & 0.289 & 0.155 & 0.731 & 0.199 & 0.851 & 0.098 & 0.905 & 0.004 & 0.905 & 0.047 & 0.502 & 0.256 & 0.905 & 0.004 & 0.906 & 0.001 & 0.898 & 0.000 \\ 
 -3,-1 & 0.175 & 0.152 & 0.533 & 0.229 & 0.572 & 0.160 & 0.684 & 0.021 & 0.672 & 0.084 & 0.363 & 0.225 & 0.682 & 0.024 & 0.689 & 0.005 & 0.674 & 0.000 \\ 
 -1/2,-2 & 0.012 & 0.010 & 0.062 & 0.084 & 0.226 & 0.060 & 0.875 & 0.009 & 0.299 & 0.174 & 0.035 & 0.022 & 0.218 & 0.106 & 0.741 & 0.115 & 0.681 & 0.095 \\ 
 -1,-2 & 0.015 & 0.011 & 0.019 & 0.022 & 0.032 & 0.023 & 0.144 & 0.100 & 0.078 & 0.092 & 0.015 & 0.023 & 0.039 & 0.077 & 0.258 & 0.204 & 0.121 & 0.126 \\ 
 -3,-2 & 0.015 & 0.012 & 0.017 & 0.025 & 0.031 & 0.031 & 0.098 & 0.097 & 0.056 & 0.087 & 0.015 & 0.019 & 0.043 & 0.054 & 0.182 & 0.158 & 0.188 & 0.201 \\ 

 \hline
\end{tabu}$$
}
\end{table}
\end{landscape}

\begin{landscape}
\begin{table}
\caption{Scaled MISE values and sd values where $n=2^{12}$}{\fontsize{6pt}{6pt}\selectfont
$$\begin{tabu}[c]{|c||cccccc|cccccc|cccccc|}
   \hline
    {\rm Pareto} & \multicolumn{6}{c|}{m=n^{1/4}} & \multicolumn{6}{c|}{m=n^{1/2}} & \multicolumn{6}{c|}{m=n^{3/4}} \\
 \ell & {\rm PE} & {\rm sd} & {\rm CV} & {\rm sd}  & {\rm PI} & {\rm sd}  & {\rm PE} & {\rm sd} & {\rm CV} & {\rm sd} & {\rm PI} & {\rm sd} & {\rm PE} & {\rm sd} & {\rm CV} & {\rm sd} & {\rm PI} & {\rm sd}   \\ \hline
1/2 & 0.107 & 0.451 & 10.05 & 0.083 & 10.34 & 15.39 & 0.433 & 0.804 & 34.50 & 51.04 & 5735 & 173000 & 64.82 & 1061 & 17500 & 44480 & 10.83 & 0.000 \\ 
  1 & 0.004 & 0.004 & 1.682 & 0.043 & 3.820 & 3.328 & 0.042 & 0.070 & 5.073 & 23.29 & 208.8 & 842.2 & 8.390 & 141.0 & 49.05 & 206.5 & 11880 & 187200 \\ 
  3 & 0.002 & 0.002 & 0.034 & 0.009 & 0.043 & 0.041 & 0.017 & 0.020 & 0.515 & 0.179 & 0.653 & 0.621 & 0.964 & 5.703 & 8.932 & 6.583 & 9.987 & 12.90 \\ 
  10 & 0.002 & 0.002 & 0.034 & 0.007 & 0.019 & 0.005 & 0.013 & 0.016 & 0.311 & 0.073 & 0.183 & 0.048 & 0.414 & 1.135 & 3.073 & 1.143 & 1.782 & 0.582 \\ 
  
\hline\hline
 {\rm T} dist. & \multicolumn{6}{c|}{m=n^{1/4}} & \multicolumn{6}{c|}{m=n^{1/2}} & \multicolumn{6}{c|}{m=n^{3/4}} \\
\ell & {\rm PE} & {\rm sd} & {\rm CV} & {\rm sd}  & {\rm PI} & {\rm sd}  & {\rm PE} & {\rm sd} & {\rm CV} & {\rm sd} & {\rm PI} & {\rm sd} & {\rm PE} & {\rm sd} & {\rm CV} & {\rm sd} & {\rm PI} & {\rm sd}   \\ \hline
1/2 & 0.153 & 0.314 & 9.952 & 0.020 & 9.941 & 12.32 & 0.328 & 0.633 & 10.74 & 0.000 & 791.8 & 24500 & 44.45 & 758.3 & 10.83 & 0.000 & 10.83 & 0.000 \\ 
  1 & 0.015 & 0.012 & 1.474 & 0.159 & 2.010 & 1.050 & 0.038 & 0.048 & 15.07 & 60.21 & 124.9 & 425.5 & 3.296 & 20.33 & 224.5 & 1409 & 2846 & 23870 \\ 
  3 & 0.003 & 0.002 & 0.015 & 0.009 & 0.005 & 0.005 & 0.015 & 0.017 & 0.181 & 0.100 & 0.064 & 0.073 & 3.085 & 71.44 & 2.735 & 2.013 & 0.945 & 0.927 \\ 
  10 & 0.002 & 0.002 & 0.003 & 0.002 & 0.003 & 0.002 & 0.013 & 0.014 & 0.023 & 0.017 & 0.019 & 0.014 & 3.835 & 70.94 & 0.207 & 0.186 & 0.156 & 0.127 \\ 
  
 \hline\hline
 {\rm Burr} & \multicolumn{6}{c|}{m=n^{1/4}} & \multicolumn{6}{c|}{m=n^{1/2}} & \multicolumn{6}{c|}{m=n^{3/4}} \\
 c,\ell & {\rm PE} & {\rm sd} & {\rm CV} & {\rm sd}  & {\rm PI} & {\rm sd}  & {\rm PE} & {\rm sd} & {\rm CV} & {\rm sd} & {\rm PI} & {\rm sd} & {\rm PE} & {\rm sd} & {\rm CV} & {\rm sd} & {\rm PI} & {\rm sd}   \\ \hline
1/2,1/2 & 457.0 & 250.5 & 845.4 & 0.001 & 790.4 & 24.09 & 549.6 & 257.4 & 1037 & 0.000 & 1037 & 0.000 & 1183 & 3083 & 161600 & 281100 & 1070 & 0.000 \\ 
  1,1/2 & 0.053 & 0.239 & 10.06 & 0.058 & 10.81 & 12.13 & 0.467 & 0.751 & 10.74 & 0.000 & 10.74 & 0.000 & 21.49 & 73.41 & 287.8 & 2015 & 10.83 & 0.000 \\ 
  3,1/2 & 0.002 & 0.003 & 0.126 & 0.027 & 0.480 & 0.426 & 0.027 & 0.027 & 0.532 & 0.262 & 12.34 & 15.00 & 0.657 & 0.744 & 19.65 & 45.41 & 361.4 & 1623 \\ 
  1/2,1 & 0.094 & 0.362 & 11.71 & 0.060 & 11.11 & 0.180 & 0.292 & 0.437 & 11.55 & 2.249 & 10.95 & 0.000 & 56.17 & 327.3 & 5921 & 32310 & 10.86 & 0.000 \\ 
  1,1 & 0.004 & 0.005 & 1.702 & 0.033 & 3.707 & 2.544 & 0.033 & 0.044 & 1.912 & 0.216 & 316.2 & 1121 & 28.81 & 242.3 & 33.34 & 66.91 & 22.28 & 204.4 \\ 
  3,1 & 0.002 & 0.002 & 0.008 & 0.004 & 0.009 & 0.007 & 0.016 & 0.016 & 0.124 & 0.043 & 0.131 & 0.119 & 1.263 & 8.059 & 1.842 & 0.618 & 1.983 & 0.715 \\ 
  1/2,3 & 0.009 & 0.007 & 0.460 & 0.576 & 4.466 & 3.297 & 0.022 & 0.024 & 2.125 & 15.44 & 360.0 & 1408 & 0.948 & 1.459 & 21.50 & 30.83 & 3822 & 34410 \\ 
  1,3 & 0.002 & 0.003 & 0.033 & 0.009 & 0.042 & 0.036 & 0.017 & 0.018 & 0.500 & 0.137 & 0.630 & 0.420 & 0.633 & 1.195 & 8.755 & 7.714 & 9.593 & 8.013 \\ 
  3,3 & 0.002 & 0.001 & 0.003 & 0.002 & 0.003 & 0.002 & 0.012 & 0.012 & 0.027 & 0.021 & 0.024 & 0.016 & 0.459 & 0.749 & 0.253 & 0.170 & 0.227 & 0.121 \\ 
  
   \hline\hline
  {\rm Fr\acute{e}chet} & \multicolumn{6}{c|}{m=n^{1/4}} & \multicolumn{6}{c|}{m=n^{1/2}} & \multicolumn{6}{c|}{m=n^{3/4}} \\
 \gamma & {\rm PE} & {\rm sd} & {\rm CV} & {\rm sd}  & {\rm PI} & {\rm sd}  & {\rm PE} & {\rm sd} & {\rm CV} & {\rm sd} & {\rm PI} & {\rm sd} & {\rm PE} & {\rm sd} & {\rm CV} & {\rm sd} & {\rm PI} & {\rm sd}   \\ \hline
  5 & 11140 & 3586 & 13700 & 0.000 & 12880 & 631.1 & 10390 & 3020 & 13700 & 0.000 & 13700 & 0.000 & 31820 & 236400 & 698200 & 1287000 & 13700 & 0.000 \\ 
  2 & 0.108 & 0.511 & 10.80 & 0.075 & 12.45 & 66.10 & 0.389 & 0.636 & 407.6 & 447.7 & 10.84 & 0.060 & 71.17 & 1018 & 50850 & 96120 & 10.84 & 0.000 \\ 
  1 & 0.004 & 0.004 & 1.654 & 0.088 & 4.228 & 3.256 & 0.042 & 0.057 & 5.166 & 29.94 & 290.8 & 1724 & 4.465 & 32.78 & 138.6 & 1138 & 16710 & 307700 \\  
  1/2 & 0.002 & 0.002 & 0.027 & 0.010 & 0.168 & 0.221 & 0.018 & 0.019 & 0.528 & 0.227 & 3.440 & 6.561 & 3.172 & 34.75 & 15.59 & 22.17 & 57.71 & 240.0 \\ 
1/4 & 0.002 & 0.002 & 0.022 & 0.013 & 0.009 & 0.004 & 0.015 & 0.017 & 0.253 & 0.144 & 0.117 & 0.035 & 0.676 & 4.942 & 3.613 & 2.658 & 1.590 & 0.540 \\  
      
  \hline\hline
    {\rm Weibull} & \multicolumn{6}{c|}{m=n^{1/4}} & \multicolumn{6}{c|}{m=n^{1/2}} & \multicolumn{6}{c|}{m=n^{3/4}} \\
 \kappa & {\rm PE} & {\rm sd} & {\rm CV} & {\rm sd}  & {\rm PI} & {\rm sd}  & {\rm PE} & {\rm sd} & {\rm CV} & {\rm sd} & {\rm PI} & {\rm sd} & {\rm PE} & {\rm sd} & {\rm CV} & {\rm sd} & {\rm PI} & {\rm sd}   \\ \hline
1/2 & 0.007 & 0.004 & 0.039 & 0.010 & 0.303 & 0.131 & 0.015 & 0.016 & 0.526 & 0.148 & 4.235 & 3.119 & 0.756 & 6.083 & 6.578 & 3.672 & 47.28 & 97.62 \\ 
  1 & 0.002 & 0.002 & 0.031 & 0.007 & 0.015 & 0.005 & 0.012 & 0.013 & 0.242 & 0.060 & 0.115 & 0.035 & 0.433 & 1.330 & 1.932 & 0.586 & 0.924 & 0.326 \\ 
  3 & 0.002 & 0.001 & 0.002 & 0.002 & 0.002 & 0.002 & 0.011 & 0.012 & 0.014 & 0.010 & 0.014 & 0.011 & 0.490 & 3.103 & 0.092 & 0.073 & 0.090 & 0.081 \\ 
  10 & 0.002 & 0.002 & 0.002 & 0.002 & 0.002 & 0.002 & 0.011 & 0.013 & 0.019 & 0.014 & 0.018 & 0.012 & 0.465 & 2.997 & 0.126 & 0.108 & 0.134 & 0.101 \\ 
  
   \hline\hline
 {\rm rev. Burr} & \multicolumn{6}{c|}{m=n^{1/4}} & \multicolumn{6}{c|}{m=n^{1/2}} & \multicolumn{6}{c|}{m=n^{3/4}} \\
 c,\ell & {\rm PE} & {\rm sd} & {\rm CV} & {\rm sd}  & {\rm PI} & {\rm sd}  & {\rm PE} & {\rm sd} & {\rm CV} & {\rm sd} & {\rm PI} & {\rm sd} & {\rm PE} & {\rm sd} & {\rm CV} & {\rm sd} & {\rm PI} & {\rm sd}   \\ \hline
 -1/2,-1/3 & 33250 & 107.8 & 33300 & 0.000 & 33300 & 0.000 & 54800 & 0.000 & 54800 & 0.000 & 54800 & 0.000 & 59610 & 35.55 & 59600 & 0.000 & 59600 & 0.000 \\ 
  -1,-1/3 & 37.26 & 0.079 & 37.40 & 0.000 & 37.40 & 0.000 & 23.44 & 4.776 & 26.50 & 0.000 & 26.50 & 0.000 & 22.32 & 3.069 & 25.40 & 0.000 & 25.40 & 0.000 \\ 
 -3,-1/3 & 2.457 & 0.161 & 2.670 & 0.000 & 2.719 & 0.003 & 0.141 & 0.258 & 1.300 & 0.000 & 1.300 & 0.000 & 0.423 & 1.253 & 1.040 & 0.000 & 1.040 & 0.000 \\ 
 -1/2,-1 & 2.056 & 0.706 & 2.849 & 0.002 & 1.873 & 0.498 & 1.173 & 0.992 & 3.260 & 0.000 & 3.260 & 0.000 & 1.597 & 0.843 & 3.340 & 0.000 & 3.340 & 0.000 \\ 
  -1,-1 & 0.226 & 0.246 & 0.493 & 0.272 & 0.493 & 0.272 & 0.023 & 0.035 & 0.904 & 0.006 & 0.904 & 0.006 & 0.467 & 0.845 & 0.898 & 0.000 & 0.898 & 0.000 \\ 
 -3,-1 & 0.076 & 0.109 & 0.364 & 0.220 & 0.364 & 0.220 & 0.009 & 0.010 & 0.684 & 0.016 & 0.684 & 0.016 & 0.498 & 1.580 & 0.674 & 0.000 & 0.674 & 0.000 \\ 
 -1/2,-2 & 0.015 & 0.061 & 0.013 & 0.008 & 0.013 & 0.008 & 0.017 & 0.022 & 0.221 & 0.048 & 0.221 & 0.048 & 0.379 & 0.441 & 0.876 & 0.009 & 0.876 & 0.009 \\ 
 -1,-2 & 0.002 & 0.001 & 0.012 & 0.016 & 0.012 & 0.016 & 0.012 & 0.011 & 0.028 & 0.035 & 0.028 & 0.035 & 0.383 & 0.449 & 0.095 & 0.091 & 0.095 & 0.091 \\ 
 -3,-2 & 0.006 & 0.003 & 0.012 & 0.019 & 0.012 & 0.019 & 0.011 & 0.011 & 0.035 & 0.050 & 0.035 & 0.050 & 0.298 & 0.370 & 0.165 & 0.186 & 0.165 & 0.186 \\ 

 \hline
\end{tabu}$$
}
\end{table}
\end{landscape}

\end{document}